\documentclass[12pt]{article}
\usepackage[final]{epsfig}
\usepackage{graphics}
\usepackage{amsmath}
\usepackage{amsfonts}
\usepackage{latexsym}
\usepackage{amssymb}
\usepackage{graphicx}
\usepackage{epstopdf}
\usepackage{multicol,multirow}
\newtheorem{lemma}{Lemma}[section]
\newtheorem{proposition}[lemma]{Proposition}
\newtheorem{remark}[lemma]{Remark}

\newtheorem{theorem}[lemma]{Theorem}
\newtheorem{definition}[lemma]{Definition}

\newtheorem{proposition-conjecture}[lemma]{Proposition-conjecture}
\newtheorem{problem}[lemma]{Problem}

\begin{document}
\newcommand{\eps}{{\varepsilon}}
\newcommand{\proofend}{\hfill$\Box$\bigskip}
\newcommand{\C}{{\mathbb C}}
\newcommand{\Q}{{\mathbb Q}}
\newcommand{\R}{{\mathbb R}}
\newcommand{\Z}{{\mathbb Z}}
\newcommand{\RP}{{\mathbb {RP}}}
\newcommand{\CP}{{\mathbb {CP}}}
\newcommand{\PP}{{\mathbb {P}}}
\newcommand{\ep}{\epsilon}
\newcommand{\G}{{\Gamma}}

\def\proof{\paragraph{Proof.}}
\def\prooflem{\paragraph{Proof of lemma.}}
\def\proofthm{\paragraph{Proof of theorem.}}
\def\proofprop{\paragraph{Proof of proposition.}}


\newcommand{\marginnote}[1]
{
}

\newcounter{bk}
\newcommand{\bk}[1]
{\stepcounter{bk}$^{\bf BK\thebk}$%
\footnotetext{\hspace{-3.7mm}$^{\blacksquare\!\blacksquare}$
{\bf BK\thebk:~}#1}}

\newcounter{fs}
\newcommand{\fs}[1]
{\stepcounter{fs}$^{\bf FS\thefs}$%
\footnotetext{\hspace{-3.7mm}$^{\blacksquare\!\blacksquare}$
{\bf FS\thefs:~}#1}}


\title {The Pentagram map in higher dimensions \\and KdV flows}

\author{Boris Khesin\thanks{
School of Mathematics, Institute for Advanced Study, Princeton, NJ 08540, USA and Department of Mathematics,
University of Toronto, Toronto, ON M5S 2E4, Canada;
e-mail: \tt{khesin@math.toronto.edu}
}
\,  and Fedor Soloviev\thanks{
Department of Mathematics,
University of Toronto, Toronto, ON M5S 2E4, Canada;
e-mail: \tt{soloviev@math.toronto.edu}
}
\\
}

\date{}

\maketitle

\begin{abstract}
We extend the definition of the pentagram map from 2D to higher dimensions and describe its integrability properties for both closed and twisted polygons by presenting its Lax form.
The corresponding continuous limit of the pentagram map in dimension $d$ is shown to be the
$(2,d+1)$-equation of the KdV hierarchy, generalizing the Boussinesq equation in 2D.
\end{abstract}


\section*{Introduction} \label{intro}

The pentagram map was originally defined  in \cite{Schwartz} as a map on plane convex polygons defined up to their projective equivalence, where a new polygon is spanned by the shortest diagonals of the initial one, see Fig.~1. It was shown to exhibit quasi-periodic behaviour under iterations. This map was extended to the case of twisted polygons and its integrability in 2D was proved in  \cite{OST99}. The integrability of the pentagram map in the case of closed polygons was proved in  \cite{FS11, OST11}.

\begin{figure}[hbtp!]\label{2d-hexagon}
\centering
\includegraphics[width=1.8in]{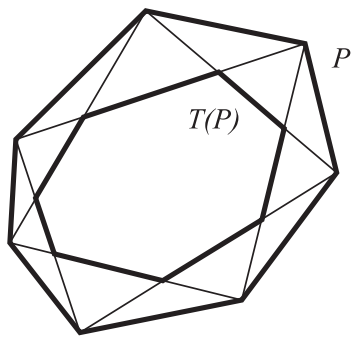}
\caption{\small The image $T(P)$ of a heptagon $P$ ($n=7; \, M=\mathrm{Id}$) under the 2D  pentagram map.}
\label{leaves}
\end{figure}

In this paper we extend this definition of the 2D pentagram map to any dimension 
and describe the 3D case in more detail.
We prove algebro-geometric integrability of the pentagram map by presenting it as a discrete zero-curvature equation (a Lax-type equation, which
implies Arnold--Liouville integrability), and study the continuous limit of the map. 
We refer to \cite{KS} for more results, proofs, and details of the constructions described in this announcement.

Note that a higher dimensional generalization for the class of corrugated polygons, which have the property that their consecutive diagonals are not skew but do intersect,  was treated in \cite{GSTV}.   
 A variety of possible higher-dimensional generalizations with an integrable continuous limit was considered in \cite{Beffa}. 
In this paper we propose a definition of ``diagonal hyperplanes" for generic  higher dimensional polygons, which leads to integrable systems not only in the continuous limit but to genuine
discrete integrable systems and can be regarded as natural integrable discretizations
of the KdV-type equations.
We start with a treatment of the 3D case and describe the $d$-dimensional case later.
We outline the geometry of the pentagram map in the real setting, but complexify the spaces and maps to describe the algebro-geometric integrability. 

\bigskip


\section{Integrability of the 3D pentagram map}

We start by extending the set of space polygons to include so-called twisted ones, whose ends are related by a monodromy operator:

\begin{definition}
{\rm
A {\it twisted $n$-gon} in a projective space $\RP^3$ with a monodromy $M \in SL_4$
is a map $\phi: \Z \to \RP^3$, such that
$\phi(k+n) =  M \circ \phi(k)$ for each $k\in \Z$ and where $M$ acts naturally on $\RP^3$.
Two twisted $n$-gons are {\it equivalent} if there is a transformation $g \in SL_4$ such that $g \circ \phi_1=\phi_2$.
}
\end{definition}
We assume that the vertices $\phi(k), \; k \in \Z,$ are in general position, i.e., in particular, no $4$ consecutive vertices
of an $n$-gon belong to one and the same $2$-dimensional plane in $\RP^3$. The following pentagram map $T$ is generically defined on the space ${\mathcal P}_n$
of twisted $n$-gons considered up to the above equivalence:

\begin{definition}
{\rm
Given an $n$-gon $\phi$ in  $\RP^3$, for each $k\in \Z$ consider the two-dimensional
``short-diagonal plane" $P_k:=(\phi(k-2), \phi(k), \phi(k+2))$ passing through 3 vertices $\phi(k-2), \phi(k), \phi(k+2)$.
Take the intersection point of the three consecutive planes $P_{k-1}, P_k, P_{k+1}$ and call it the image of
the vertex $\phi(k)$ under the {\it space pentagram map} $T$. (We assume the general position, so that  every
three consecutive planes $P_k$ for the given $n$-gon intersect at a point.)
}
\end{definition}
\begin{figure}[hbtp!]
\centering
\includegraphics[width=2.5in]{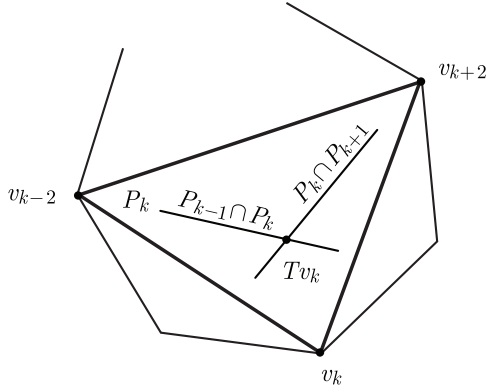}
\caption{\small The image $Tv_k$ of the vertex $v_k=\phi(k)$ in $\PP^3$.}
\end{figure}

For simplicity, in 3D we discuss only the case of odd $n$ in this research announcement.
(Even values of $n$ require more delicate treatment, for which we refer the reader to \cite{KS}.)
The coordinates on the space ${\mathcal P}_n$ are introduced in the following way.

It turns out that for odd $n$ there exists the unique lift  of the vertices $\phi(k) \in \RP^3,\; k \in \Z,$
of a given $n$-gon 
to the vectors $V_k \in \R^4,\; k \in \Z,$ satisfying the conditions of nondegeneracy
$\det(V_j, V_{j+1}, V_{j+2}, V_{j+3})=1$ and quasi-periodicity $V_{j+n}=MV_j$ for all $ j \in \Z,$ where
$M\in SL_4$ is the monodromy matrix. 
These vectors satisfy the difference equations
$$
V_{j+4} = a_j V_{j+3} + b_j V_{j+2} + c_j V_{j+1} - V_j,\; j \in \Z,
$$
where the sequences $a_j,b_j,c_j$ are $n$-periodic.
The numbers $(a_j,b_j,c_j),\; 0 \le j \le n-1,$ form a system of  $3n$ coordinates on the space of twisted $n$-gons ${\mathcal P}_n$. One can also introduce  ``local"
coordinates on ${\mathcal P}_n$ based on cross-ratios   and somewhat similar to the so-called $(x,y)$-coordinates in the 2D case, see \cite{KS}.

\medskip

\subsection{A discrete zero-curvature equation}

Algebro-geometric integrability of the pentagram map is based on a zero-curvature equation with a spectral parameter.
In the discrete case, it is an equation of the form
\begin{equation}\label{zero-eq}
 L_{i,t+1}(\lambda) = P_{i+1,t}(\lambda) L_{i,t}(\lambda) P_{i,t}^{-1}(\lambda),
\end{equation}
which represents a dynamical system.
This equation may be regarded as a compatibility condition of an over-determined system of equations:
 \begin{equation*}
 \begin{cases}
  L_{i,t}(\lambda) \Psi_{i,t}(\lambda) = \Psi_{i+1,t}(\lambda)\\
  P_{i,t}(\lambda) \Psi_{i,t}(\lambda) = \Psi_{i,t+1}(\lambda),
 \end{cases}
 \end{equation*}
for an auxiliary function  $\Psi_{i,t}(\lambda)$.
To describe this equation we complexify the pentagram map and consider it over $\C$.

The above has an analogue in the continuous case: 
a zero curvature equation $\partial_t L-\partial_x P=[P,L]$
is a compatibility condition which provides the existence of an auxiliary function $\psi=\psi(t,x)$ 
satisfying a system  of differential equations $\partial_x \psi= L\psi$ and $\partial_t \psi= P\psi\,.$

\begin{theorem}\label{lax-abc}
The 3D pentagram map on twisted $n$-gons with  odd $n$ admits a zero-curvature representation with the Lax function $L_{j,t}(\lambda)$ given by
 \[
 L_{j,t}(\lambda) =
 \begin{pmatrix}
  c_j/\lambda & 1/\lambda & 0 & 0\\
  b_j         & 0        & 1  & 0\\
  a_j/\lambda & 0        & 0  & 1/\lambda\\
  -1          & 0        & 0  & 0
 \end{pmatrix} =
 \begin{pmatrix}
  0       & 0   & 0       & -1\\
  \lambda & 0   & 0       & c_j\\
  0       & 1   & 0       & b_j\\
  0       & 0   & \lambda & a_j
 \end{pmatrix}^{-1}
 \]
in the coordinates $a_j,b_j,c_j,\,0\le j\le n-1$, and an appropriate matrix function  $P_{j,t}(\lambda)$ 
satisfying (\ref{zero-eq}), where $\lambda \in \C$ is the spectral parameter.
\end{theorem}

In a sense, equation (\ref{zero-eq}) implies all our integrability results in the discrete case, while 
the explicit form of $P_{i,t}$ is not important for the present exposition. 

\medskip

The pivotal property responsible for integrability of the pentagram map is its property of scaling invariance.
In the 3D case in the coordinates $(a_j,b_j,c_j),\; 0 \le j \le n-1,$ this means that 
the pentagram map $T$ is invariant with respect to the transformations
$a_j \to a_j s, \; b_j \to b_j, \; c_j \to c_j s,\; 0 \le j \le n-1.$  The invariance follows from the explicit formulas of the map.
We prove that the existence of equation (\ref{zero-eq}) and the corresponding formula for $L_{j,t}(\lambda)$
follow from such a scaling invariance. The spectral parameter $\lambda$ is related to the scaling parameter $s$ via $\lambda = 1/s^2.$

\subsection{Spectral curve}

In the discrete case, there exist  analogues of monodromy and Floquet--Bloch solutions:
\begin{definition}
 {\rm
 {\it The monodromy operators} $M_{i,t},\, i=0,...,n-1,$ are defined as
 the following ordered products of the corresponding Lax functions:
$
M_{i,t} = L_{i+n-1,t} L_{i+n-2,t} ... L_{i+1,t} L_{i,t},
$
where the (integer) index $t$  represents the time variable.
}
\end{definition}

\begin{definition}\label{ijg-def}
{\rm
For an odd $n$ define the {\it spectral function} $R(\lambda,k)$ as
$$
R(\lambda,k):=\det{(M_{0,0}(\lambda) - k I)}.
$$
The {\it spectral curve} $\Gamma$ is the normalization of the compactification of the curve $R(\lambda,k)=0$.
}
\end{definition}

Note that one could define the spectral function with help of any of the monodromies:
$R(\lambda,k):=\det{(M_{i,t}(\lambda) - k I)}$. Indeed, by definition, the monodromy operators for different $i$ are conjugate to each other, while 
 equation (\ref{zero-eq}) implies that the monodromy operators for different $t$ are also conjugate (i.e., satisfy a discrete Lax equation):
$$
M_{i,t+1}(\lambda) = P_{i,t}(\lambda) M_{i,t}(\lambda) P_{i,t}^{-1}(\lambda).
$$
This implies that the spectral function $R(\lambda,k)$
is  an  invariant of the pentagram map.
After a multiplication by a suitable power of $\lambda$, the equation $R(\lambda,k)=0$ becomes a polynomial relation between $\lambda$ and $k$,
whereas its coefficients are integrals of motion for the pentagram map.
Namely, we define the {\it integrals of motion} $I_j,J_j,G_j,\; 0 \le j \le q=[n/2],$ as
the coefficients of the expansion
$$
R(\lambda,k) = k^4 - k^3 \left( \sum_{j=0}^q G_j \lambda^{j-n} \right) + k^2 \left( \sum_{j=0}^q J_j \lambda^{j-q-n} \right)
- k \left( \sum_{j=0}^q I_j \lambda^{j-2n} \right) + \lambda^{-2n} = 0.
$$
The  set  $R(\lambda,k)=0$
is an algebraic curve in $\C^2$. A standard procedure (adding the infinite points and normalization with
a few blow-ups) makes it into a compact Riemann surface,  the corresponding  spectral curve 
$\Gamma$. For a generic  spectral curve the integrals of motion  are independent  
  polynomials in the $(a,b,c)$-coordinates.

\begin{theorem}\label{thm:genus}
For an odd $n$ the genus $g$ of the spectral curve $\Gamma$ is $g = 3q$, where $q=[n/2]$.
\end{theorem}

The spectral curve and its Jacobian, a natural torus  associated with it, are  starting points for the algebraic-geometric integrability.  It turns out that one can recover
Lax functions from the spectral curve and a point on the Jacobian, and vice versa: this correspondence is locally one-to-one.

\subsection{Pentagram dynamics in spectral data}

\begin{definition}
{\rm
A {\it Floquet--Bloch solution} $\psi_{i,t}$ of a difference equation $\psi_{i+1,t} = L_{i,t} \psi_{i,t}$ is an eigenvector
of the monodromy operator:
$$
M_{i,t} \psi_{i,t} = k \psi_{i,t}.
$$
To make them uniquely defined, we normalize the vectors $\psi_{i,t}, \, t \ge 0$ so that the sum of their components  is
equal to $1$ and denote the normalized vectors  by $\Bar{\psi}_{i,t}$, i.e., $\Bar{\psi}_{i,t} := \psi_{i,t}/\left(\sum_{j=1}^4 \psi_{i,t,j}\right)$.
}
\end{definition}

\begin{theorem}\label{bloch-th}
A Floquet--Bloch solution $\Bar{\psi}_{i,t}$ is a meromorphic vector function on $\Gamma$.
Its pole divisor $D_{i,t}$ has degree $g+3$.
\end{theorem}

\begin{definition}
{\rm
Let $J(\Gamma)$ be the Jacobian of the spectral curve $\Gamma$, and $[D_{i,t}]$
be the equivalence class of the divisor $D_{i,t}$ under the Abel map.
}
\end{definition}

\begin{theorem}\label{spectral-th}
 For an odd $n$, the spectral map
 $S: \mathcal{P}_n = (a_i,b_i,c_i, 0 \le i \le n-1) \to (\Gamma, [D_{i,t}])$ is non-degenerate at a generic point, i.e., locally, it is one-to-one.
 
 The pentagram dynamics corresponds to the motion of the pole divisor $[D_{i,t}]$ along the corresponding Jacobian $J(\Gamma)$.
The equivalence class $[D_{i,t}]\in J(\Gamma)$ has the following linear time evolution along the Jacobian:
 $$
 [D_{i,t}] = [D_{0,0} - t(O_1+O_3) + i(O_2+O_3) + (t-i)(W_1+W_2)] \in J(\Gamma),
 $$
 where the points $W_1,W_2 \in \Gamma$ correspond to $\lambda = \infty$, and
 $O_1,O_2,O_3 \in \Gamma$ are the points above $\lambda = 0$ (the point $O_1$ corresponds to a finite $k$, whereas $O_2$ and $O_3$ correspond to infinite $k$, and 
   $O_2$ is a branch point).
 \end{theorem}

This theorem implies that for an odd $n$ the time evolution in $J(\Gamma)$ happens 
to be along a straight line by a constant shift. This way the above theorem  describes 
the time evolution of the pentagram map and proves its algebraic-geometric integrability.

For an even $n$ (which we deal with in \cite{KS}) the dynamics turns out to be 
more complicated: the evolution  goes along a ``staircase'', i.e., its square is a constant shift. 
This dichotomy is similar to the 2D case, see \cite{FS11}.


\medskip

\subsection{Closed polygons}

 Closed polygons in $\PP^3$  correspond to the monodromy $M=\pm\mathrm{Id} \in SL_4$ and
 form a subspace of codimension $15=\dim SL_4$ in the space of all twisted polygons $\mathcal{P}_n$ of dimension $3n$. Such a monodromy   corresponds
to the spectral curves with either $(\lambda,k)=(1,1)$ or $(\lambda,k)=(1,-1)$ being a quadruple point.
More precisely, we have:

\begin{theorem}\label{thm:closed}
 Closed polygons  in $\PP^3$ are singled out by the condition that either $(\lambda,k)=(1,1)$
 or $(\lambda,k)=(1,-1)$ is a quadruple point of \, $\Gamma$.
 The genus of \, $\Gamma$ drops to $g=3q-6$ for $q=[n/2]$ and odd $n$.
 The dimension of the Jacobian $J(\Gamma)$ drops by $6$ for closed polygons.
 Theorem~\ref{bloch-th} holds with this genus adjustment, and Theorem~\ref{spectral-th} holds
 verbatim for closed polygons. 
 \end{theorem}
 
 \begin{remark}
{\rm 
 The algebraic conditions implying that $(1,\pm 1)$ is a quadruple point are:
 \begin{itemize}
 \item $R(1,\pm 1) = 0$,
 \item $\partial_k R(1,\pm 1) = \partial_\lambda R(1,\pm 1) = 0$,
 \item $\partial_k^2 R(1,\pm 1) = \partial_\lambda^2 R(1,\pm 1) = \partial_{k \lambda}^2 R(1,\pm 1) = 0$,
 \item $\partial_k^3 R(1,\pm 1) = \partial_\lambda^3 R(1,\pm 1) = \partial_{kk \lambda}^3 R(1,\pm 1)= \partial_{k \lambda \lambda}^3 R(1,\pm 1) = 0$.
 \end{itemize}
 One can show that among these 10 linear equations on the integrals of motion there are
 only 9 independent ones due to the following relation:
 $$
 R(1,\pm 1)=\pm\partial_k R(1,\pm 1) - \dfrac{1}{2} \partial_k^2 R(1,\pm 1) \pm \dfrac{1}{6} \partial_k^3 R(1,\pm 1).
 $$
 At the same time, the dimension of the Jacobian of the spectral curve drops by $6$. 
 Thus the subspace ${\mathcal C}_n$ has codimension $15$, which matches the above calculation of dimensions.
 }
 \end{remark}

\subsection{An invariant symplectic structure}\label{S:ham-struc}

To describe an invariant symplectic structure on leaves of the space 
of twisted polygons ${\mathcal P}_n$ in 3D we employ the Krichever--Phong universal formula
\cite{KP97,KP98}. In the 2D case such a symplectic structure was shown to coincide with that 
induced on these leaves by the invariant Poisson structure found in \cite{OST99}, see \cite{FS11}. 
The description below is somewhat implicit and, in a sense, universal: it is applicable in the higher-dimensional cases of $\PP^d$ with $d>3$ as well. 
Finding an explicit expression, e.g., in the coordinates $(a_i,b_i,c_i)$,
of the symplectic structure or of the corresponding Poisson structure is still an open problem.

 \begin{definition}[\cite{KP97,KP98}]\label{KP-universal}
{\rm
 Krichever--Phong's universal formula defines a {\it pre-symplectic form} on the space of Lax operators, i.e.,
 on the space $\mathcal{P}_n$.
 It is given by the expression:
 $$
 \omega = -\dfrac{1}{2} \sum_{\lambda=0,\infty} {\text{res}} \thinspace \text{Tr}\left( \Psi_0^{-1} M_0^{-1} \delta M_0 \wedge \delta \Psi_0
 \right) \dfrac{d\lambda}{\lambda}.
 $$
 The matrix $\Psi_{0}(\lambda)$ is composed of the normalized eigenvectors $\Bar{\psi}_{0,0}$
  on different sheets of $\Gamma$ over the $\lambda$-plane,
 and it diagonalizes the matrix $M_{0}=M_{0,0}$. (In this definition we drop the second index, since
 all variables correspond to the same moment $t$.)

 The {\it leaves} of the 2-form $\omega$ are defined as submanifolds of $\mathcal{P}_n$, where the expression $\delta \ln{k}\, (d\lambda/\lambda)$
 is holomorphic. The latter expression is considered as a 1-form on the spectral curve $\Gamma$.
}
 \end{definition}

 \begin{theorem}\label{symp-leaves}\label{thm:rank}
 For an odd $n$ the leaves of the 2-form $\omega$ in $\mathcal{P}_n$ 
 are singled out by the $3$ conditions
 $$
 \delta G_0 = \delta I_0 = \delta J_q = 0,
 $$
 where $I_0,G_0,J_q$ are the integrals of motion defined above.
The restriction of $\omega$ to these leaves is well-defined (i.e., independent of the normalization of 
the Floquet--Bloch solutions $\Psi_0$) and non-degenerate, and hence symplectic. 
 This symplectic form is invariant with respect to the pentagram map, that is the 
 evolution given by the Lax equation.

 The rank of the invariant 2-form $\omega$  restricted to these leaves 
  is equal to $2g$, where $g$ is the genus of the spectral curve $\Gamma$, i.e.,  $g = 3[n/2]$.
 \end{theorem}
 
Recall that the dimension of the space $\mathcal{P}_n$ is $6q+3$, where $q=[n/2]$.
Since the codimension of the leaves is $3$, their dimension matches the doubled  dimension of the tori: $2g=6q$. We also note that the Arnold--Liouville integrability in the complex case 
implies integrability in the real one, since the formula for the symplectic structure is algebraic.


\section{The pentagram maps in higher dimensions}

\subsection{Definition for any dimension}

First we extend the notion of a twisted $n$-gon to an arbitrary dimension $d$.

\begin{definition}
{\rm
A {\it twisted $n$-gon} in a projective space $\RP^d$ with a monodromy $M \in SL_{d+1}$
is a map $\phi: \Z \to \RP^d$, such that
$\phi(k+n) =  M \circ \phi(k)$ for each $k\in \Z$ and where $M$ acts naturally on $\RP^d$.
Two twisted $n$-gons are {\it equivalent} if there is a transformation $g \in SL_{d+1}$ such that $g \circ \phi_1=\phi_2$.

Similarly to the 3D case, we assume that the vertices $v_k:=\phi(k), \; k \in \Z,$ are in general position (i.e., no $d+1$ consecutive vertices lie in the same hyperplane in $\RP^d$), and denote by   ${\mathcal P}_n$ the space of twisted $n$-gons considered up to the above equivalence. 

For a generic twisted $n$-gon in $\RP^d$ define the ``short-diagonal" hyperplane  $P_k$ passing
through $d$ vertices of the $n$-gon by taking $d$-tuple consisting every other vertex and centered 
at a given vertex $v_k$. Namely, for  odd dimension $d=2\varkappa+1$
we consider the {\it short-diagonal hyperplane} $P_k$ through the $d$ vertices
$$
P_k:=(v_{k-2\varkappa}, v_{k-2\varkappa+2},...,v_k, ..., v_{k+2\varkappa}),
$$
while for even dimension $d=2\varkappa$ we take $P_k$  passing through the  $d$ vertices
$$
P_k:=(v_{k-2\varkappa+1}, v_{k-2\varkappa+3}, ..., v_{k-1},v_{k+1}, ..., v_{k+2\varkappa-1}).
$$
}
\end{definition}

The following pentagram map $T$ is generically defined on the space ${\mathcal P}_n$ of twisted $n$-gons:

\begin{definition}
{\rm
The {\it higher pentagram map} $T$ takes a vertex $v_k$ of a generic twisted
$n$-gon  in  $\PP^d$
to the intersection point of the $d$ consecutive short-diagonal planes $P_i$ around $v_k$.
Namely, for odd $d=2\varkappa+1$ one takes the intersection of the planes
$$
Tv_k=P_{k-\varkappa}\cap P_{k-\varkappa+1}\cap ...\cap P_k\cap ...\cap P_{k+\varkappa}\,,
$$
while for even $d=2\varkappa$ one takes
the intersection of the planes
$$
Tv_k=P_{k-\varkappa+1}\cap P_{k-\varkappa+2}\cap ...\cap P_{k}\cap ...\cap P_{k+\varkappa}\,.
$$
As usual, we assume that the vertices are in ``general position,'' and every $d$ 
consecutive hyperplanes $P_i$ intersect at one point in $\PP^d$.
The  map $T$ is well defined on the equivalence classes of  generic $n$-gons in $\PP^d$.
}
\end{definition}

\begin{remark}\label{diff-eq}
{\rm
One can show that there exists a unique lift of the vertices $v_k=\phi(k) \in \RP^d$
to the vectors $V_k \in \R^{d+1}$ satisfying $\det|V_j, V_{j+1}, ..., V_{j+d}|=1$ and $ V_{j+n}=MV_j,\; j \in \Z,$ where
$M\in SL_{d+1}$, if and only if the condition $gcd(n,d+1)=1$ holds.
The corresponding difference equations have the form
\begin{equation}\label{eq:difference_anyD}
V_{j+d+1} = a_{j,d} V_{j+d} + a_{j,d-1} V_{j+d-1} +...+ a_{j,1} V_{j+1} +(-1)^{d} V_j,\quad j \in \Z,
\end{equation}
with $n$-periodic coefficients in the index $j$.
It allows one to introduce coordinates $\{ a_{j,k} ,\;0\le j\le n-1, \; 1\le k\le d \}$ on the space of twisted $n$-gons in $\RP^d$.
}
\end{remark}


\subsection{Complete integrability}

It turns out that the pentagram map defined this way has a special scaling invariance,
which implies the existence of a Lax representation, and allows one to prove its integrability.
First, we complexify the spaces of $n$-gons and the map.

\begin{proposition-conjecture} {\bf (The scaling invariance)}
The pentagram map on twisted $n$-gons in $\CP^d$ is invariant with respect to the following scaling transformations:
\begin{itemize}
\item for odd $d=2\varkappa+1$ the transformations are
$$
a_{j,1} \to s a_{j,1},\; a_{j,3} \to s a_{j,3},\;a_{j,5} \to s a_{j,5},\;...\;, a_{j,d} \to s a_{j,d}\,,
$$
while other coefficients $a_{j,2l}$ with $l=1,...,\varkappa$ do not change;

\item for even $d=2\varkappa$ the transformations are
$$
a_{j,1} \to s^{-1}a_{j,1} ,\; a_{j,2} \to s^{-2}a_{j,2} ,\;...\:, a_{j,\varkappa} \to  s^{-\varkappa}a_{j,\varkappa},
$$
$$
a_{j,\varkappa+1} \to s^\varkappa a_{j,\varkappa+1} ,\; a_{j,\varkappa+2} \to s^{\varkappa-1}a_{j,\varkappa+2}, \;...\:, a_{j,d-1} \to s^2 a_{j,d-1} ,
\; a_{j,d} \to s a_{j,d}
$$
\end{itemize}
for all $s\in \C$.
\end{proposition-conjecture}

We proved this proposition up to dimension $d\le 6$ by studying the explicit formulas for the pentagram map, but have no general proof for  $d > 6$. It would be very interesting to find it.

\begin{problem}
Find a general proof of the scaling invariance of the pentagram map in any dimension $d$.
\end{problem}

For the following theorem we assume this conjecture on scaling invariance.

\begin{theorem}\label{thm:lax_anyD}
The scale-invariant pentagram map on twisted $n$-gons in any dimension $d$ is a completely integrable system. It is described by the  Lax matrix
\[
L_{j,t}(\lambda) =
\left(
\begin{array}{cccc|c}
0 & 0 & \cdots & 0    &(-1)^d\\ \cline{1-5}
\multicolumn{4}{c|}{\multirow{4}*{$D(\lambda)$}} & a_{j,1}\\
&&&& a_{j,2}\\
&&&& \cdots\\
&&&& a_{j,d}\\
\end{array}
\right)^{-1},
\] 
where $D(\lambda)$ is the following diagonal $(d \times d)$-matrix:
\begin{itemize}
\item for odd $d=2\varkappa+1$, one has
$D(\lambda) = \text{diag}(\lambda,1,\lambda,1,...,\lambda)$;

\item for even $d=2\varkappa$, one has
\[
D_{ii}(\lambda)=
\begin{cases}
1,       &\text{if $i \ne \varkappa+1$,}\\
\lambda, &\text{if $i = \varkappa+1$.}
\end{cases}
\]
\end{itemize}
\end{theorem}


{\it Sketch of proof.} Rewrite the difference equation  \eqref{eq:difference_anyD} in the matrix form.
It is described by the transformation matrix
\[
N_j :=
\left(
\begin{array}{ccc|c}
0 &  \cdots & 0    &(-1)^d\\ \cline{1-4}
\multicolumn{3}{c|}{\multirow{3}*{\rm{Id}}} & a_{j,1}\\
&&& \cdots\\
&&& a_{j,d}\\
\end{array}
\right),
\] 
where $\mathrm{Id}$ is the identity $(d\times d)$-matrix. 
Then the monodromy  $M$ for  twisted $n$-gons is the product $M=N_0 N_1...N_{n-1}$.
Note that the pentagram map defined 
on classes of projective equivalence preserves the conjugacy class of $M$.
E.g., assume that $d$ is odd. Then using the scale invariance, replace 
$a_{j,2k+1}$ by $s a_{j,2k+1}$ for all $k$ in the right column to obtain new matrices $N_j(s)$. The  pentagram map preserves the conjugacy class of  the new monodromy $M(s):=N_0(s) ...N_{n-1}(s)$ for any $s$. 
Although $N_j(s)$ could have already been taken as a Lax matrix $L_j(s)$, for computations it is convenient to set 
$L_{j,t}^{-1}(\lambda) := \left(g^{-1} N_j(s) g \right)/s$ for  $g = \text{diag}(1,s,1,s,...,1,s)$ and $\lambda := s^{-2}$. \hfill$\Box$


\subsection{Continuous limit of the pentagram map}

Consider the continuous limit of polygons and the pentagram map on them.
In the limit  $n\to\infty$ a generic twisted $n$-gon becomes a smooth non-degenerate quasi-periodic curve $\gamma(x)$ in $\RP^d$.
Its lift $G(x)$ to $\R^{d+1}$ is defined by the conditions that the components of the vector function
$G(x)=(G_1,...,G_{d+1})(x)$ provide the homogeneous coordinates for 
$\gamma(x)=(G_1:...:G_{d+1})(x)$ in $\RP^d$ and $\det|G(x),G'(x),...,G^{(d)}(x)|=1$ for all $x\in \R$.
Furthermore, $G(x+2\pi)=MG(x)$ for a given $M\in SL_{d+1}$. Then $G(x)$
satisfies the linear differential equation of order $d+1$:
$$
G^{(d+1)}+u_{d-1}(x)G^{(d-1)}+...+u_1(x)G'+u_0(x)G=0
$$
with periodic coefficients $u_i(x)$, which is a continuous limit of difference equation \eqref{eq:difference_anyD}. (Here $'$ stands for $d/dx$.)

Fix a small $\epsilon>0$ and consider the case of odd $d=2\varkappa+1$. 
A continuous analog of the hyperplane $P_k$ is the hyperplane $P_\epsilon(x)$ passing through
$d$ points $\gamma(x-\varkappa\epsilon),..., \gamma(x),...,\gamma(x+\varkappa\epsilon)$ of the curve $\gamma$.

Let $\ell_\epsilon (x)$ be the envelope curve for the family of hyperplanes $P_\epsilon(x)$  for a fixed $\epsilon$.
The envelope condition means that  $P_\epsilon(x)$ are the osculating hyperplanes of the curve $\ell_\epsilon (x)$, that is  the point $\ell_\epsilon (x)$ belongs to the plane $P_\epsilon(x)$, while the vector-derivatives $\ell'_\epsilon (x),...,\ell^{(d-1)}_\epsilon (x)$ span this plane
 for each $x$. It means that the lift of $\ell_\epsilon (x)$ to $L_\epsilon (x)$ in $\R^{d+1}$
satisfies the system of $d=2\varkappa+1$ equations (see Fig.~3 for $d=3$):
$$
\det | G(x-\varkappa\epsilon), G(x-(\varkappa-1)\ep ), ... , G(x),..., G(x+\varkappa\epsilon),  L^{(j)}_\epsilon(x) |=0,\quad j=0,...,d-1.
$$

\begin{figure}[hbtp]
\centering
\includegraphics[width=3.5in]{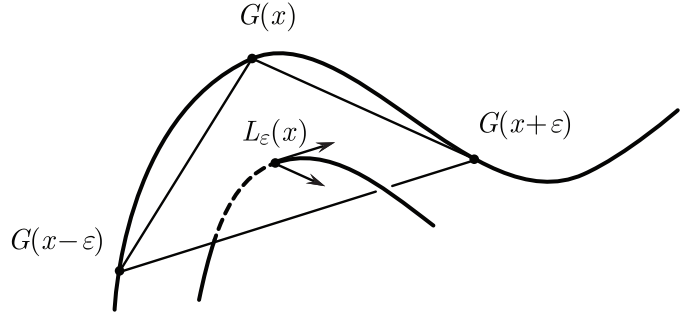}
\caption{\small The envelope $L_\ep(x)$ in 3D. The point $L_\ep(x)$ and the vectors $L_\ep'(x)$ and $L_\ep''(x)$ belong to the plane $(G(x), G(x+\ep), G(x-\ep))$.}
\end{figure}

Similarly, for even $d=2\varkappa$ the lift $L_\epsilon (x)$ satisfies the system of $d$ equations:
\begin{align*} 
\det | G(x-(2\varkappa-1)\epsilon), &\, G(x-(2\varkappa-3)\ep ),  ... , G(x-\ep),G(x+\ep),...\\
 ...,
& \, G(x+(2\varkappa-1)\epsilon),  L^{(j)}_\epsilon(x) |=0,\quad j=0,...,d-1.
\end{align*}

A continuous limit of the pentagram map $T$ is defined as the evolution of the curve $\gamma$ in the direction of the envelope $\ell_\epsilon$, as $\epsilon$ changes.
Namely, the expansion of $L_\epsilon(x)$ has the form
$$
L_\epsilon(x)=G(x)+\epsilon^2 B(x)+{\mathcal O} (\epsilon^4)
$$
and satisfies the family of differential equations:
$$
L_\epsilon^{(d+1)}+u_{d-1,\ep}(x)L_\epsilon^{(d-1)}+...+u_{1,\ep}(x)L_\epsilon'+u_{0,\ep}(x)L_\epsilon=0, \text{ where } u_{j,0}(x)=u_{j}(x).
$$
Then the corresponding expansion of the coefficients $u_{j,\ep}(x)$ as $u_{j,\ep}(x)=u_{j}(x)+\ep^2w_j(x)+{\mathcal O}(\ep^4)$,
defines the continuous limit of the pentagram map as the system of evolution differential equations $du_j(x)/dt\, =w_j(x)$ for $j=0,...,d-1$.

\medskip

\begin{theorem}
The continuous limit of the pentagram map $T$ in dimension $d$ defined by the system $du_j(x)/dt\, =w_j(x), \, j=0,...,d-1$ for $x\in S^1$  is the $(2, d+1)$-KdV flow of the Adler-Gelfand-Dickey  hierarchy on the circle.
\end{theorem}

Recall that the $(m, d+1)$-KdV flow is defined on  linear differential operators  
 $L= \partial^{d+1} + u_{d-1}(x) \partial^{d-1} + u_{d-2}(x) \partial^{d-2} + ...+ u_1(x) \partial + u_0(x)$  of order $d+1$ with periodic coefficients $u_j(x)$, where $\partial^{k}$ stands for $d^k/dx^k$.
One can define its fractional power
$L^{m/{d+1}}$ as a pseudo-differential operator for any positive integer $m$ and take
its pure differential part  $Q_m :=(L^{m/{d+1}})_+$. In particular, for $m=2$ one has $Q_2= \partial^2 + \dfrac{2}{d+1}u_{d-1}(x) $. Then the $(m, d+1)$-KdV equation is the  evolution equation on (the coefficients of) $L$ given by $dL/dt= [Q_m,L] .$

\begin{remark}
{\rm
 For $d=2$ the  (2,3)-KdV equation is the classical Boussinesq equation, found in \cite{OST99}.
Apparently, the $(2, d+1)$-KdV equation is a very robust continuous limit.
One obtains it not only for the pentagram map defined by taking every other vertex,
but also for a non-symmetric choice of vertices for the plane $P_k$.
Also, the same limit was obtained in \cite{Beffa} for a map defined by taking intersections of various planes, rather than the envelopes.
}
\end{remark}

The above scaling has a clear meaning in the continuous limit (in 2D this was proved in \cite{OST99}):

\begin{proposition}\label{cts-scaling}
The continuous limit of the scaling transformations
corresponds to the spectral shift $L\rightarrow L+\lambda$ of the differential operator $L$.
\end{proposition}

\medskip

{\bf Acknowledgments}.
We are grateful to M.~Gekhtman and S.~Tabachnikov for useful discussions and to Olga Solovieva for help with drawing the figures.
B.K. was partially supported by the Simonyi Fund and an NSERC research grant.


\end{document}